\documentclass[a4paper, 12pt]{article}

\usepackage[koi8-r,cp1251]{inputenc}
\usepackage{amsfonts,amssymb,amsmath, hyperref}
\usepackage[final]{epsfig}
\usepackage{graphicx}

\begin{document}

\begin{center}
\textbf{ERRATUM\\
to\\
 "Addendum to "Boundedness of Hausdorff operators on Hardy spaces   over locally compact groups" \\
 (J. Math. Anal. Appl. 473(2019) 519-533)"}
\end{center}

\begin{center}
A. R. Mirotin
\end{center}

\begin{center}
amirotin@yandex.ru
\end{center}

\

Abstract. We give a correct reformulation of Lemma 2 and the main result of our note "Addendum to "Boundedness of Hausdorff operators on Hardy spaces   over locally
compact groups" [J. Math. Anal. Appl. 473(2019) 519-533]", J. Math. Anal. Appl. 2019; 479(1):
872-874 (see below).

\

In the next  "Addendum" Lemma 2 and Theorem  should be reformulated as follows.

\textbf{Lemma 2}. \textit{There is a left invariant metric $\rho_1$  which is compatible with the
topology of $G$ such that every automorphism $A\in \mathrm{Aut}(G)$  is Lipschitz with respect to
every left invariant metric $\rho$  that is strongly equivalent to  $\rho_1$. Moreover, one can
choose the Lipschitz constant to be
$$
L_A=\kappa_\rho \mathrm{mod}A
$$
where the constant $\kappa_\rho$  depends on the metric $\rho$ only and  $\kappa_{\rho_1}=1$.
}

The proof of Lemma 2 above is exactly the same as for Lemma 2 in   "Addendum" below.

\

Lemma 2  implies that the following condition

\textit{for every automorphism $A(u),$  for every $x\in G,$ and for every  $r>0$ there exist a positive number $k(u)$ which depends on $u$ only and a point $x'=x'(x,u,r)\in G$ such that}
$$
A(u)^{-1}(B(x,r))\subseteq B(x',k(u)r) \eqno(*)
$$
from \cite{Hausd} holds for the metric space  $(G,\rho)$ with $x'=A(u)^{-1}(x)$ and $k(u)= \kappa_\rho /\mathrm{mod}(A(u)).$

Thus  the main result of  "Addendum" below should be  as follows.

\textbf{Theorem.} \textit{Let a left invariant metric $\rho$  be as in Lemma 2,  the doubling condition holds for the
corresponding metric measure space  $(G,\rho,\nu)$, and $k(u)=\kappa_\rho
/\mathrm{mod}(A(u)).$ For $\Phi\in L^1(\Omega,k^sd\mu)$   the
Hausdorff operator $\mathcal{H}_{\Phi,A}$ is bounded on the real
Hardy space $H^1(G)$ and
$$
\|\mathcal{H}_{\Phi,A}\|\leq C_\nu \|\Phi\|_{L^1(\Omega,k^sd\mu)}.
$$
}

Due to the Lemma 2 above the proof of this theorem is exactly the same as for main result in  \cite{Hausd}.

\begin{center}
{\Large Addendum to "Boundedness of  Hausdorff  operators
on real Hardy spaces $H^1$ over locally compact groups" }
\end{center}

In the author's paper \cite{Hausd} (see below), results by  Liflyand and collaborators on the boundedness of Hausdorff operators on the Hardy space $H^1$
over finite-dimensional real space are generalized to the case of  locally compact groups that are spaces of homogeneous type. In the following
$G$ stands for a locally compact  group which is a space of homogeneous type with respect to some quasi-metric $\rho$  and left Haar measure $\nu$.
Let $B(x,r)$ denote a quasi-ball  with respect to $\rho$ centered at $x$ of radius $r$, $(\Omega,\mu)$ be some measure space, $\Phi\in L^1_{loc}(\Omega),$
and let $A(u)$ stand for a $\mu$-measurable family of automorphisms of $G$.

 The main theorem of \cite{Hausd} asserts that the Hausdorff operator
$$
(\mathcal{H}_{\Phi, A}f)(x) =\int_\Omega \Phi(u)f(A(u)(x))d\mu(u)
$$
is bounded in the Hardy space $H^1(G)$ provided that the following condition holds:

\noindent {\it For every automorphism $A(u)$  of $G,$   for every $x\in G,$ and for every  $r>0$, there exist a positive number $k(u)$ which depends on
$A(u)$ only and a point $x'=x'(x,u,r)\in G$ such that}
$$
A(u)^{-1}(B(x,r))\subseteq B(x',k(u)r). \eqno(*)
$$

The aim of this note is twofold: to prove that condition $(*)$ holds automatically, and in the case  of locally compact metrizable group
to find the best possible value of the  constant $k(u).$

Recall that for each $x\in G, t\geq 1$ and $r > 0$
$$
\nu(B(x,tr))\leq C_\nu t^{s} \nu(B(x,r)) \eqno(D)
$$
where $C_\nu$ denotes the doubling constant, $s=\log_2C_\nu$ (see, e.g., \cite[p. 76]{HK}).

\textbf{Lemma A.} \textit{Let $G$ be as above. Then $G$ is a space of homogeneous type with respect to  the measure
$\nu$ and any left invariant metric which  is  comparable with topology of $G.$}

Proof. By \cite{MiMiMiMo13} (see also \cite[Theorem 2.1]{AMi}), there exist a metric $d$ on $G$ and positive constants $a, b,$ and $\beta$ such that
$$
ad^{1/\beta}(x,y)\leq \rho(x,y)\leq bd^{1/\beta}(x,y)\eqno(1)
$$
for all $x, y\in G.$  One can assume that $d$ is left invariant (see, e.g., \cite[Theorem 2.8.3]{HiR}). By formula (1), for every $r>0$,
$$
B_d(x,r)\subseteq B(x,br^{1/\beta})\subseteq B_d(x,\left(\frac{b}{a}\right)^\beta r) \eqno(2)
$$
(here $B_d$ denotes a ball with respect to $d$). Now, using (2) and (D) (with $t=2^{1/\beta}b/a$), we get, for every $R>0$,
\begin{align*}
\nu(B_d(x,2R))&\leq \nu(B(x,b(2R)^{1/\beta}))=\nu(B(x,(\frac{b}{a}2^{1/\beta})(aR^{1/\beta})))\\ &\leq
C_\nu(\frac{b}{a}2^{1/\beta})^s\nu(B(x,b((\frac{a}{b})^\beta R)^{1/\beta}))\leq
C_\nu(2^{1/\beta}\frac{b}{a})^s\nu(B_d(x,R)).
\end{align*}
Thus, the doubling condition holds for $d$.  Since every   metric $q$,  which is compatible with the topology of $G$, is equivalent to $d,$
the same arguments as above show that the doubling condition holds for $q$ as well.

Without loss of generality, we shall assume in the following that $\rho$ \textit{is a left invariant metric}, which is compatible with the topology of $G.$

\textbf{Lemma B.} \textit{Every automorphism  $A\in \mathrm{Aut}(G)$ is Lipschitz. Moreover, one can choose the Lipschitz constant to be}
$$
L_A=\kappa_\rho \mathrm{mod}(A),
$$
\textit{where the constant $\kappa_\rho$ depends on the metric $\rho$ only.
}

Proof. By \cite[Lemma 1]{Str},  there is a decreasing family $V_n$ of open neighborhoods of identity in $G$ such that the left invariant metric
$$
\rho_1(x,y):=\sup_n\nu(xV_n\bigtriangleup yV_n)
$$
is compatible with the topology of $G.$

We have for every $A\in \mathrm{Aut}(G),$  $x,y\in G$,
\begin{align*}
\rho_1(A(x),A(y))&=\sup_n\nu(A(x)V_n\bigtriangleup A(y)V_n)=\sup_n\nu(A(xV_n\bigtriangleup yV_n))\\
&=\mathrm{mod}(A)\sup_n\nu(xV_n\bigtriangleup yV_n)=\mathrm{mod}(A)\rho_1(x,y).
\end{align*}
Since every two  metrics which are compatible with the topology of $G$ are equivalent, the result follows. Indeed, if
$$
C_1\rho_1(x,y)\leq \rho(x,y)\leq C_2\rho_1(x,y),
$$
then the Lipschitz constant $L_A=(C_2/C_1)  \mathrm{mod}(A).$ This
value is the best possible, as the case $\rho=\rho_1$ shows.

\textbf{Remark.} There is a  very simple proof of the Lipschitzness of any automorphism $A$ of $G.$ (Indeed,
the formula $p(x,y):= \rho(A(x),A(y))$ defines a left invariant metric in $G$ which induces the original topology and thus $p$ and $\rho$ are equivalent.)
But we are interested in the \textit{best possible} Lipschitz constant, and Lemma 2 makes the job. Moreover, it gives us an explicit formula for this constant.

Lemma B  implies that condition $(\ast)$ holds with $x'=A(u)^{-1}(x)$ and $k(u)= \kappa_\rho /\mathrm{mod}(A(u)).$
This refines the main result of \cite{Hausd} as follows.

\textbf{Theorem.}  \textit{Let $\rho$ be a left invariant metric which is compatible with the topology of $G$ and $k(u)=\kappa_\rho /\mathrm{mod}(A(u)).$
For $\Phi\in L^1(\Omega,k^sd\mu)$   the Hausdorff operator $\mathcal{H}_{\Phi,A}$ is bounded on the real Hardy space $H^1(G)$ and
$$
\|\mathcal{H}_{\Phi,A}\|\leq C_\nu \|\Phi\|_{L^1(\Omega,k^sd\mu)}.
$$
}

\

\begin{center}
\bf{{\Large Boundedness of  Hausdorff  operators
on Hardy spaces $H^1$ over locally compact groups }}\\
\vspace{5mm}
\end{center}

Abstract. {\small Results of  Liflyand and collaborators on the boundedness of Hausdorff operators on the Hardy space $H^1$ over finite-dimensional real space generalized to the case of  locally compact groups that are spaces of homogeneous type. Special cases and examples  of compact  Lie groups, homogeneous groups (in particular the Heisenberg group) and finite-dimensional  spaces over division rings are considered.}
\vspace{5mm}

\vspace{5mm}

\vspace{5mm}

\section{Introduction}

 Hausdorff operators  originated from some classical summation methods. This class of operators contains some important examples, such as Hardy operator, adjoint
Hardy operator,  the Ces\`{a}ro operator. As mentioned in \cite{RF} the Riemann-Liouville fractional integral and the
Hardy-Littlewood-Polya operator   can also be reduced to  the Hausdorff operator.  As was noted in \cite{CDFZ}  the Hausdorff operator is closely related to a Calder\'{o}n-Zygmund convolution operator, too.

  The study of general Hausdorff operators on Hardy spaces $H^1$ over  the real line was pioneered by Liflyand  and M\'{o}ricz \cite{LM}. After publication of this paper   Hausdorff
operators have attracted much attention. The multidimensional case  was considered by Lerner and Liflyand  \cite{LL}, and Liflyand  \cite{ActaSz} (the case of the space $L^p(\mathbb{R}^n)$ was studied earlier in \cite{BM}).  Hausdorff operators on  spaces $H^p(\mathbb{R})$ for $p\in (0,1)$ where considered by Kanjin  \cite{Kan}, and Liflyand and Miyachi \cite{LM}.  The survey article by E.~Liflyand  \cite{Ls} contains
main results on Hausdorff operators in various settings and bibliography up to 2013. See also \cite{CFW},  and \cite{Mor}, \cite{CFL}, \cite{CFLR}, \cite{RFW}, \cite{WF}.  The recent paper by Ruan and Fan \cite{RF} contains in particular  several sharp conditions for boundedness of Hausdorff operators on the space $H^1(\mathbb{R}^n).$

The aim of this work is to generalize results  on the boundedness of Hausdorff operators on Hardy spaces over $\mathbb{R}^n$ to the case of general locally compact groups. The main task is a distillation of results about  Hausdorff operators depending only on the group and (quasi-)metric structures. So, we consider locally compact groups that are spaces of homogeneous type in the sense of Coifman and Weiss  \cite{CW}. Special cases and examples of compact  Lie groups, homogeneous groups (in particular the Heisenberg group) and finite-dimensional  spaces over division rings are also considered.

Recall that according to \cite{CW} \textit{a space of homogeneous type} is
a quasi-metric space $\Omega$ endowed with a Borel measure $\mu$ and a quasi-metric $\rho.$
And the basic assumption relating the measure
and the quasi-metric is the existence of a constant $C$ such that
$$
\mu(B(x,2r))\leq C \mu(B(x,r))
$$
for each $x\in \Omega$ and $r > 0$ ("the doubling condition"). Here $B(x,r)$ denotes a quasi-ball of radius $r$ around $x.$ The \textit{doubling constant }  is the smallest constant
$C\geq  1$ for which the last inequality holds. We denote this constant by $C_\mu.$ Then  for each $x\in \Omega, k\geq 1$ and $r > 0$
$$
\mu(B(x,kr))\leq C_\mu k^{s} \mu(B(x,r)),\eqno(D)
$$
where $s=\log_2C_\mu$ (see, e.g., \cite[p. 76]{HK}). The number $s$  sometimes takes
the role of a “dimension” for a doubling quasi-metric measure space.

Recall also the definition of the real Hardy space $H^1(\Omega)$ associated with a space of homogeneous type $\Omega$ \cite{CW}.

First note that a function
$a$ on $\Omega$ is an ($(1, \infty)$-)\textit{atom} if

(i) the support of $a$ is contained in a ball $B(x,r);$

(ii) $\|a\|_\infty\leq\frac{1}{\mu(B(x,r))};$

(iii) $\int\limits_\Omega a(x)d\mu(x) = 0.$

By definition, the \textit{real Hardy space} $H^1(\Omega)$ consists of those functions
admitting an atomic decomposition
$$
f=\sum\limits_{j}\lambda_ja_j
$$
where the $a_j$ are atoms, and $\sum_{j}|\lambda_j|<\infty$  \cite[p. 593]{CW}.

The infimum of the
numbers  $\sum_{j}|\lambda_j|$  taken over all such representations of $f$ will be denoted by
the symbol $\|f\|_{H^1}.$

\textbf{Remark 1.} Real Hardy spaces over compact connected (not necessary quasi-metric) Abelian groups  were  defined in \cite{Indag}.

\section{The general case}

In the following $G$ stands for a locally compact $\sigma$-compact  group which is a space of homogeneous type  (in particular, a homogeneous group \cite{FS}, \cite{FR}) with respect to quasi-metric $\rho$ and left Haar measure $\mu,$  and  $A:G\to \mathrm{Aut}(G)$ a $\mu$-measurable map.

 Let
$$
L_A(G)=\{\Phi:G\to \mathbb{C}: \|\Phi\|_{L_A}:=\int\limits_G |\Phi(u)|\mathrm{mod}(A(u)^{-1})d\mu(u)<\infty\}
$$
where $\mathrm{mod}(A(u)^{-1})(=1/\mathrm{mod}(A(u)))$ denotes the modulus of the automorphism $A(u)^{-1}$ (recall that the modulus of the automorphism $\varphi\in \mathrm{Aut}(G)$ satisfies $\mu(\varphi(E))=(\mathrm{mod}(\varphi))\mu(E)$ for every Borel  $E\subset G$ with finite measure, see, e.g., \cite[Chapter VII]{Burb}).

\textbf{Definition 1} (cf. \cite{BM}). Let $\Phi$ be a locally integrable function on $G.$ We define the \textit{Hausdorff  operator} with the kernel $\Phi$  by
$$
(\mathcal{H}f)(x) = (\mathcal{H}_{\Phi, A}f)(x) =\int\limits_G \Phi(u)f(A(u)(x))d\mu(u). \eqno(1)
$$

\textbf{Remark 2.}  One can assume that $A$ is defined almost everywhere on the support of $\Phi$ only.

%2) If
%$$
%\Phi(u)=\chi_{\{\mathrm{mod}(A(u))\geq 1\}}(u)\mathrm{mod}(A(u)),
%$$
%it is natural to call $\mathcal{H}_{\Phi, A}$ the Ces\`{a}ro operator on $G$ (cf. \cite[p. 127]{Ls}).

%\textbf{Definition 2.} We define the \textit{Ces\`{a}ro operator} on $G$ as follows
%$$
%(\mathcal{C}_Af)(x)=\int\limits_{\{\mathrm{mod}(A(u))\geq 1\}} f(A(u)(x))\mathrm{mod}(A(u))d\mu(u).
%$$

We need four lemmas to prove our main result.

\textbf{Lemma 1.} \textit{Let $\Phi\in L_A(G).$ Then the operator $\mathcal{H}_{\Phi, A}$ is bounded in $L^1(G)$  and}
$$
\|\mathcal{H}_{\Phi, A}\|\leq \|\Phi\|_{L_A}.
$$

Proof. Using Fubini Theorem and \cite[VII.1.4, formula (31)]{Burb} we have for $f\in L^1(G)$
$$
\|\mathcal{H}_{\Phi, A}f\|_{L^1}=\int\limits_G \left|\int\limits_G\Phi(u)f(A(u)(x))d\mu(u) \right|d\mu(x)\leq
$$
$$
\int\limits_G|\Phi(u)|\int\limits_G|f(A(u)(x))|d\mu(x)d\mu(u)=
$$
$$
\int\limits_G|\Phi(u)|\mod(A(u)^{-1})\left(\int\limits_G|f(x)|d\mu(x)\right)d\mu(u)=\|\Phi\|_{L_A}\|f\|_{L^1}.
$$

\textbf{Lemma 2.} \textit{Let $(\Omega,\rho)$ be  $\sigma$-compact quasi-metric space with positive Radon measure $\mu,$ and let $\mathcal{F}(\Omega)$ be some Banach  space of  $\mu$-measurable functions on $\Omega.$ Assume  that the convergence of a sequence strongly in  $\mathcal{F}(\Omega)$ yields the convergence of some subsequence to the same function for $\mu$-almost all $x\in \Omega.$ Let $F(u,x)$ be a function such that $F(u,\cdot)\in \mathcal{F}(\Omega)$ for $\mu$-almost all $u\in \Omega$ and $u\mapsto F(u,\cdot):\Omega\to \mathcal{F}(\Omega)$ is Bochner integrable with respect to $\mu.$ Then for $\mu$-almost all} $x\in \Omega$
$$
\left((B)\int\limits_\Omega F(u,\cdot)d\mu(u)\right)(x)=\int\limits_\Omega F(u,x)d\mu(u).
$$

Proof. Let $K_m$ be an increasing sequence of compact subsets of $\Omega$ and $\Omega=\cup_{m=1}^\infty K_m.$ Then $\mu(K_m)<\infty$ and
$$
(B)\int\limits_\Omega F(u,\cdot)d\mu(u)=\lim\limits_{m\to\infty}(B)\int\limits_{K_m}F(u,\cdot)d\mu(u).
$$
By \cite{MiMiMiMo13} (see also \cite[Theorem 2.1]{AMi}) there exist a metric $d$ on $\Omega$ and positive constants $a, b,$ and $\beta$ such that
$$
ad^{1/\beta}(x,y)\leq \rho(x,y)\leq bd^{1/\beta}(x,y)\eqno(3)
$$
for all $x, y\in \Omega.$ Since $(\Omega,\rho)$ and $(\Omega,d)$ are isomorphic as uniform spaces, Theorem 1 from \cite{vN} remains true for $(\Omega,\rho)$ along with its proof. Therefore \cite[p. 203]{vN} for every $m$ there are a sequence of partitions $P^{(n)}=(\Omega_j^{(n)})_{j=1}^{N(n)}$ of $K_m$ with  the property $\max_j\mathrm{diam}(\Omega_j^{(n)})\to 0$ as $n\to\infty$ and a sequence of sample point sets
$S^{(n)}=\{u^{(n)}_j:j=1,2,\dots, N(n)\}$ such that
$$
(B)\int\limits_{K_m}F(u,\cdot)d\mu(u)=\lim\limits_{n\to\infty}\sum\limits_{j=1}^{N(n)}F(u^{(n)}_j,\cdot)\mu(\Omega_j^{(n)})
$$
strongly in $\mathcal{F}(\Omega),$ and therefore the sequence in the right-hand side contains a subsequence that converges to the function in the left-hand side $\mu$-almost everywhere.
This implies that for $\mu$-almost all $x\in \Omega$
$$
(B)\int\limits_{K_m}F(u,\cdot)d\mu(u)(x)=\int\limits_{K_m}F(u,x)d\mu(u)
$$
and lemma 2 follows.

\textbf{Lemma 3.} \textit{There are such $a, b>0$ that for all $x, x'\in G$ and} $r>0$
$$
\mu(B(x',\frac{a}{b}r))\leq \mu(B(x,r))\leq \mu(B(x',\frac{b}{a}r)).
$$

Proof. Let positive constants $a, b, \beta,$ and a metric $d$ on $G$ be such that (3) is valid. One can assume that $d$ is left invariant (see, e.g., \cite[Theorem 2.8.3]{HiR}). By formula (3)
$$
B_d(x,\left(\frac{r}{b}\right)^\beta)\subseteq B(x,r)\subseteq B_d(x,\left(\frac{r}{a}\right)^\beta)
$$
and therefore
$$
\mu(B_d(x,\left(\frac{r}{b}\right)^\beta))\leq \mu(B(x,r))\leq \mu(B_d(x,\left(\frac{r}{a}\right)^\beta)). \eqno(4)
$$
Since $\mu$ and $d$ are left invariant, $\mu(B_d(x,R))=\mu(B_d(x',R))$ for all $x, x'\in G, R>0.$
It follows in view of (4) that
$$
\mu(B(x,r))\geq \mu(B_d(x,\left(\frac{r}{b}\right)^\beta))=\mu(B_d(x',\left(\frac{r}{b}\right)^\beta))=
$$
$$
\mu(B_d(x',\left(\frac{ar/b}{a}\right)^\beta))\geq \mu(B(x',\frac{a}{b}r)).
$$
The proof of the second inequality is similar.

Consider the following condition: for every automorphism $A(u),$  for every $x\in G,$ and for every  $r>0$ there exist a positive number $k(u)$ which depends of $u$ only and a point $x'=x'(x,u,r)\in G$ such that
$$
A(u)^{-1}(B(x,r))\subseteq B(x',k(u)r) \eqno(*)
$$
 (in fact, $k(u)$ depends of $A(u)$). In the following we choose $k$ to be a $\mu$-measurable function,  $s=\log_2C_\mu.$

We shall say that $\Phi\in L^1_{k^s}(G)$ if
$$
\|\Phi\|_{ L^1_{k^s}}:=\int\limits_{\Omega} |\Phi(u)|k(u)^sd\mu(u)<\infty.
$$

\textbf{Lemma 4.} \textit{If the  condition $(*)$ holds, then} $L^1_{k^s}(G)\subseteq L_{A}(G).$

Proof. For every $x, u\in G,$ and for every  $r>0$ the  condition $(*)$ implies that
$$
\mathrm{mod}(A(u)^{-1})\mu(B(x,r))=\mu(A(u)^{-1}(B(x,r)))\leq \mu(B(x',k(u)r)).
$$
On the other hand, we have by Lemma 3 and formula (D)
$$
\mu(B(x',k(u)r))\leq \mu(B(x,\frac{b}{a}k(u)r))\leq C_\mu k(u)^s\left(\frac{b}{a}\right)^s\mu(B(x,r)).
$$
Then $\mathrm{mod}(A(u)^{-1})\leq C_\mu (b/a)^s k(u)^s$ and the desired inclusion  follows.

Now we are in position to prove our main theorem.

\textbf{Theorem 1.} \textit{Let the  condition $(*)$ holds. For $\Phi\in L^1_{k^s}(G)$   the Hausdorff operator $\mathcal{H}$ is bounded on the real Hardy space $H^1(G)$ and}
$$
\|\mathcal{H}\|\leq C_\mu\left(\frac{b}{a}\right)^s\|\Phi\|_{ L^1_{k^s}}.
$$

Proof. We use the approach from \cite{ActaSz}. First note that by lemmas 4 and 1 the integral in (1) exists. Since for $f\in H^1(G)$ we have $\|f\|_{L^1}\leq \|f\|_{H^1},$ one can apply lemma 2 and   formula (1) can be rewritten as follows:
$$
\mathcal{H}_{\Phi, A}f =\int\limits_G \Phi(u)f\circ A(u)d\mu(u),
$$
the Bochner integral with respect to $H^1$ norm, (as usual, $\circ$ denotes the composition operation) and therefore
$$
\|\mathcal{H}_{\Phi, A}f\|_{H^1}\leq \int\limits_G |\Phi(u)|\|f\circ A(u)\|_{H^1}d\mu(u). \eqno(5)
$$
We wish to estimate the right-hand side of (5) from above by using $(*).$ If $f$ has an atomic decomposition $f=\sum_{j}\lambda_ja_j,$ then
$$
f\circ A(u)=\sum\limits_{j}\lambda_ja_j\circ A(u).\eqno(6)
$$
We claim that  $a_{j,u}':=C_\mu^{-1}(bk(u)/a)^{-s}a_j\circ A(u)$ is an atom, as well. Indeed, the condition $(*)$ shows that if $a_{j}$ is supported in $B(x_j,r_j),$ the function $a_{j,u}'$ is supported in $B(x_j',k(u)r_j),$ and thus (i) holds for $a_{j,u}'.$ Next, by lemma 3 and the doubling condition,
$$
\mu(B(x_j',k(u)r_j))\leq \mu(B(x_j,\frac{b}{a}k(u)r_j))\leq C_\mu \left(\frac{b}{a}k(u)\right)^s\mu(B(x_j,r_j)).
$$
Then
$$
\|a_j\circ A(u)\|_\infty\leq\frac{1}{\mu(B(x_j,r_j))}\leq  C_\mu \left(\frac{b}{a}k(u)\right)^s\frac{1}{\mu(B(x_j',k(u)r_j))},
$$
and (ii) is also valid  for $a_{j,u}'.$ Finally, the cancelation property (iii)  for $a_{j,u}'$ follows from \cite[VII.1.4, formula (31)]{Burb}.

Since by (6)
$$
f\circ A(u)=\sum\limits_{j}\left(C_\mu \left(\frac{b}{a}k(u)\right)^s\lambda_j\right)a_{j,u}',
$$
we get
$$
\|f\circ A(u)\|_{H^1}\leq C_\mu \left(\frac{b}{a}k(u)\right)^s\sum\limits_{j}|\lambda_j|.
$$
Therefore $\|f\circ A(u)\|_{H^1}\leq C_\mu \left(bk(u)/a\right)^s\|f\|_{H^1}$ and the conclusion  of the theorem follows from the formula (5).

\textbf{Remark 3.} If  $\rho$ is a  left invariant quasi-metric, one can take  $a=b=1$ in theorem 1 because in this case lemma 3 holds trivially with such $a$ and $b$.

\vspace{5mm}

\section{Special cases and examples}
\

\subsection{Compact Lie groups}
\

As mentioned in \cite[p. 588, Example (7)]{CW} compact
Lie groups with natural distances and Haar measures are spaces of homogeneous type. Moreover, \textit{the condition $(*)$ holds for such groups automatically} as the following lemma shows.

\textbf{Lemma 5.} \textit{Let $G$ be a compact Lie group with left invariant metric $\rho.$ Every automorphism
$A\in \mathrm{Aut}(G)$ is Lipschitz, i.e. for some constant $k>0$ and for every} $x, y\in G$
$$
\rho(A(x), A(y))\leq k \rho(x,y).
$$

Proof.  Taking into account that every  compact Lie group is smoothly isomorphic to a matrix group, one can assume that $G$ is such a group. Let an automorphism
$A\in \mathrm{Aut}(G)$ induces an automorphism $\hat{A}$ of the Lie algebra $\mathfrak{g}$ of $G$ such that $A(\exp X)=\exp (\hat{A} X)$ (see, e.g. \cite{Chev}). Consider a sufficiently small   neighborhood $U$ of unit $e\in G$ such that $\exp^{-1}$ is defined in $U.$ For $x\in U$ let $X:=\exp^{-1}(x)$ and let the norm $\|\cdot\|$ on $\mathfrak{g}$ corresponds to the metric $\rho.$ Since (infimum below is taken over all curves $\alpha\in C^1([0,1],G)$ with  $\alpha(0)=e,$ $\alpha(1)=\exp{X}$)
$$
\rho(\exp{X}, e)=\inf\limits_\alpha\int\limits_0^1\left\|\alpha'(t)\right\|dt\leq \int\limits_0^1\left\|\frac{d}{dt}\exp{(tX)}\right\|dt=\|X\|+o(\|X\|)
$$
and therefore $\rho(A(x), e)=\rho(\exp (\hat{A} X),e)=\|\hat{A}X\|+o(\|X\|),$
we have
$$
\limsup\limits_{x\to e}\frac{\rho(A(x),e)}{\rho(x, e)}= \limsup\limits_{X\to 0}\frac{\|\hat{A}X\|}{\|X\|}\leq \|\hat{A}\|.
$$
Thus the  function $\rho(A(x),e)/\rho(x, e)$ is bounded in some open  neighborhood  $V$ of unit. Since  it is also continuous on the compact set $G\setminus V,$  it is bounded, $\rho(A(x),e)/\rho(x, e)\leq k.$ To finish the proof one should substitute $y^{-1}x$ in place of $x$ in the last inequality.

Now theorem 1 yields the following corollary (see remark 3).

\textbf{Corollary 1.}  \textit{Let $G$ be a compact Lie group with left invariant metric $\rho.$ For $\Phi\in L^1_{k^s}(G)$   the Hausdorff operator $\mathcal{H}_{\Phi, A}$ is bounded on the real Hardy space $H^1(G)$ and}
$$
\|\mathcal{H}_{\Phi, A}\|\leq C_\mu\|\Phi\|_{ L^1_{k^s}}.
$$

%Since every automorphism of a compact group has modulus 1, the  Ces\`{a}ro operator in this situation has the form
%$$
%(\mathcal{C}_Af)(x)=\int\limits_{G} f(A(u)(x))d\mu(u)
%$$
%and it is bounded provided $k^s\in L^1(G)$ (in this case $\Phi(u) =1$).

\textbf{Examples.} (1) \textit{The $n$-dimensional torus $\mathbb{T}^n.$} We assume that $\mathbb{T}^n$ is equipped  with the invariant metric $\rho(x,y)= \max_{1\leq i\leq n}d(x_i,y_i)$ (here $x=(x_i),$ $y=(y_i)\in \mathbb{T}$ and $d$ is a usual metric in $\mathbb{T}$).

The one-dimensional torus  possesses only two automorphisms $z\mapsto z$ and $z\mapsto -z.$ Therefore we can take $k(u)=1$ for every $A(u)\in \mathrm{Aut}(\mathbb{T})$ and then $L^1_{k^s}(\mathbb{T})=L^1(\mathbb{T}).$ It follows that the condition $\Phi\in L^1_{k^s}(G)$ of theorem 1 is sharp in general and that bounded Hausdorff operators on $H^1(\mathbb{T})$
turns out to be very simple: $\mathcal{H}=aI+bJ$ where $If=f, Jf(z)=f(-z),$ and $a, b\in \mathbb{R}.$
%for the Ces\`{a}ro operator we should to require also that $a, b\geq 0, a+b=1.$

In the general case $n>1$ all elements of  $\mathrm{Aut}(\mathbb{T}^n)$ have the
form
$$
A(z_1,\dots,z_n)=(z_1^{m_{11}}z_2^{m_{21}}\dots z_n^{m_{n1}},\dots,z_1^{m_{1n}}z_2^{m_{2n}}\dots z_n^{m_{nn}})
$$
where the matrix $(m_{ij})$ belongs to $\mathrm{GL}(n,\mathbb{Z})$ and $\det(m_{ij})=\pm 1$ (see, e.g., \cite[(26.18)(h)]{HiR}).

Thus for every measurable map $A:\mathbb{T}^n\to \mathrm{Aut}(\mathbb{T}^n), u\mapsto (m_{ij}(u))$ the corresponding Hausdorff operator over $\mathbb{T}^n$ takes the form
$$
(\mathcal{H}_{\Phi, A}f)(z)=\int\limits_{\mathbb{T}^n}\Phi(u)f(z_1^{m_{11}(u)}\dots z_n^{m_{n1}(u)},\dots,z_1^{m_{1n}(u)}\dots z_n^{m_{nn}(u)})d\mu_n(u),
$$
where $\mu_n$ denotes the normalized Lebesgue measure on $\mathbb{T}^n.$

In this example the measure of the ball $B(x,r)$ for sufficiently small $r>0$ has the form $c_nr^n,$ and therefore $C_\mu=2^n.$

(2)  \textit{The special unitary group }$SU(2).$ It is a compact connected Lie group which  is isomorphic to the group of unit quaternions and it is known that all automorphisms of $SU(2)$ are inner. It follows that $k(u)=1$ for every $A(u)\in \mathrm{Aut}(SU(2))$ (we consider   a bi-invariant metric in $SU(2)$) and therefore $L^1_{k^s}(SU(2))=L^1(SU(2))$. Every Hausdorff   operator for $SU(2)$ has the form
 (below $b:SU(2)\to SU(2)$ is a
$\mu$-measurable map)
$$
(\mathcal{H}_{\Phi, b}f)(x)=\int\limits_{SU(2)}\Phi(u)f(b(u)xb(u)^{-1})d\mu(u).
$$
According to theorem 1 this operator is bounded in $H^1(SU(2))$ if  $\Phi\in L^1(SU(2)).$

The group $SU(2)$  as a space of homogeneous type may be identified with the 3-sphere $S^3\subset \mathbb{R}^4$ endowed with the natural distance  and volume (action of $SU(2)$ preserves the inner product in $\mathbb{C}^2$). It follows that $C_\mu=8$ in this example.

%Thus, the  Ces\`{a}ro operator in this example is bounded and has the form
%$$
%(\mathcal{C}_bf)(x)=\int\limits_{SU(2)}f(b(u)xb(u)^{-1})d\mu(u).
%$$

\

\subsection{Homogeneous groups}
\

According to \cite{FS} a \textit{homogeneous group} $G$
is a connected simply connected Lie group whose Lie algebra is equipped
with dilations. It induces the dilation structure $D_\lambda$ ($\lambda >0$) on the  group
$G$ such that $D_\lambda \in\mathrm{ Aut}(G)$ \cite[p. 5, 6]{FS} (see also \cite{FR}).

The group $G$ is endowed with   a homogeneous (quasi-)norm,  a continuous nonnegative function $|\cdot|$ on $G$
 which  satisfies  $|x^{-1}| = |x|,$
$|D_\lambda(x)| = \lambda |x|$ for all $x\in G,$  $\lambda > 0,$ and  $|x| = 0$ if and only if $x=e,$ the unit of $G.$ Moreover, the formula $\rho(x,y):=|y^{-1}x|$ defines a left invariant quasi-metric on $G$ \cite[p. 9, Proposition 1.6]{FS}.

Let $\mu$ be a (bi-invariant) Haar measure on $G$ normalized in such a way that $\mu(B(x,r))=r^Q$
where $Q$ is the so called \textit{homogeneous dimension} of $G$ \cite[p. 10]{FS}. Then the doubling condition holds and $C_{\mu}=2^Q$ (see also \cite[Lemma 3.2.12]{FR}).

Let $\lambda:G\to (0,\infty)$ be any $\mu$-measurable function. Then the family of automorphisms $A(u):=D_{\lambda(u)}$ enjoys the property $(*)$ with $k(u)=1/\lambda(u).$  Indeed, since $D_\lambda^{-1}(x^{-1})=(D_\lambda^{-1}(x))^{-1},$  we have for every $\lambda>0$
$$
D_\lambda^{-1}(B(x,r))=\{D_\lambda^{-1}(y): |yx^{-1}|<r\}=\{z: |D_\lambda(z)x^{-1}|<r\}=
$$
$$
\{z: |D_\lambda(zD_\lambda^{-1}(x^{-1}))|<r\}=
\{z: \lambda|zD_\lambda^{-1}(x^{-1})|<r\}=
$$
$$
\{z: |z(D_\lambda^{-1}(x))^{-1})|<r/\lambda\}=B(D_\lambda^{-1}(x),r/\lambda).
$$

Since $\mu(B(x,r))=r^Q,$ it follows also that $\mathrm{mod}(D_\lambda)=\lambda^Q.$

\textbf{Definition 2.} We define the Hausdorff operator $\mathcal{H}_{\Phi,\lambda}$ for the homogeneous group $G$ via the formula (1) with $A(u)=D_{\lambda(u)}.$

%It particular the corresponding Ces\`{a}ro operator takes the form (see definition 2)
%$$
%(\mathcal{C}_\lambda f)(x)=\int\limits_{\{\lambda(u)\geq 1\}}f(D_{\lambda(u)}(x))\lambda(u)^Qd\mu(u).
%$$

Then theorem 1 and remark 3 imply the following

  \textbf{Corollary 2.} \textit{The operator $\mathcal{H}_{\Phi,\lambda}$ is bounded on $H^1(G)$ provided $\Phi\in L^1_{1/\lambda^{Q}}(G)$ and}
$$
\|\mathcal{H}_{\Phi,\lambda}\|\leq 2^Q\|\Phi\|_{ L^1_{1/\lambda^{Q}}}.
$$
%\textit{It particular the Ces\`{a}ro operator $\mathcal{C}_\lambda$  is bounded on $H^1(G)$ provided}
%$\mu(\{\lambda(u)\geq 1\})<\infty$ \textit{and in this case}

%$$
%\|\mathcal{C}_\lambda\|\leq 2^Q\left(\frac{b}{a}\right)^Q\mu(\{\lambda(u)\geq 1\}).
%$$

\textbf{Examples.} (3)\textit{ Heisenberg groups} (see, e.g.,  \cite{Foll}). If $n$ is a positive integer, the Heisenberg
group $\mathbf{H}_n$ is the group whose underlying manifold is $\mathbb{R}^n\times \mathbb{R}^n\times\mathbb{R}$ and whose
multiplication is given by ($v, w, v', w'\in \mathbb{R}^n,  t,t'\in \mathbb{R}$)
$$
(v,w,t)(v',w',t') = \left(v+v', w+w', t+t'+\frac{1}{2}(v\cdot w'-w\cdot v')\right)
$$
($v\cdot w$ stands for  the
usual inner product on $\mathbb{R}^n$). Then $\mathbf{H}_n$ is a homogeneous group with dilations
$$
D_\lambda(v,w,t) = (\lambda v, \lambda w, \lambda^2 t)
$$
(there are another families of dilations on $\mathbf{H}_n,$  see  \cite[p.  7]{FS} where the isomorphic version of $\mathbf{H}_n$ is considered). The Haar measure of $\mathbf{H}_n$ is the Lebesgue measure $dudvdt$ of $\mathbb{R}^{2n+1},$ and the homogeneous dimension of $\mathbf{H}_n$ equals to $2n+2$ (see, e.g., \cite[p. 642]{Semmes}). The left invariant Heisenberg distance $d_H$ on $\mathbf{H}_n$ is derived from the  homogeneous norm $|(v,w,t)|_H:=c_n((v^2+w^2)^2+t^2)^{1/4}$ (with an appropriate constant $c_n$ which guarantee the relation $\mu(B(x,r))=r^Q$).  So, corollary 2 is valid for $\mathbf{H}_n$ with  $Q=2n+2.$

\textbf{Remarks 4.} 1) There are  automorphisms of $\mathbf{H}_n$ distinct from $D_\lambda,$  see  \cite[Chapter I, Theotem (1.22)] {Foll} for the description of all  automorphisms of $\mathbf{H}_n$.  So, one can define the Hausdorff operator for $\mathbf{H}_n$ by definition 1 using this description. Thus,  $\mathcal{H}_{\Phi,\lambda}$ is a special case of  Hausdorff operator for $\mathbf{H}_n$ in a sense of definition 1.
 Using  automorphisms of $\mathbf{H}_n$ generated by the real symplectic  group $Sp(n, \mathbb{R})$ (see  \cite[p. 20] {Foll}) we can define another  special case of Hausdorff operator for $\mathbf{H}_n$  as follows. Consider the measurable map $S:\mathbf{H}_n\to Sp(n, \mathbb{R})$. Then the  corresponding Hausdorff operator  takes the form
$$
(\mathcal{H}_{\Phi, S}f)(v,w,t) =\int\limits_{\mathbb{R}^n\times \mathbb{R}^n\times\mathbb{R}} \Phi(v',w',t')f(S(v',w',t')(v,w),t)dv'dw'dt'.
$$

2) Special cases of Hausdorff operator on $\mathbf{H}_n$ were considered in \cite{RFW}.

(4) \textit{Strict upper triangular groups}  \cite[p. 6, 7]{FS}, \cite[Subsection VII.3.3]{Burb}. Let $T_1(n,\mathbb{R})$ be the group of all $n\times n$ real
matrices $(a_{ij})$ such that $a_{ii} = 1$ for $1 \leq i \leq n$ and $a_{ij} = 0$ when $i > j.$
Then $T_1(n,\mathbb{R})$ is a homogeneous group with Haar measure $\mu(d(a_{ij}))=\otimes_{i<j}da_{ij}$ and dilations
$$
D_\lambda(a_{ij}) = \left(\lambda^{j-i}a_{ij}\right).
$$
 It is known that $\mu(dD_\lambda (a_{ij}))=\lambda^Q\mu(d(a_{ij}))$ \cite[p. 10]{FS}. Since $\otimes_{i<j}\lambda^{j-i}da_{ij}=\lambda^Q \otimes_{i<j}da_{ij}$ where
 $$
 Q=\sum\limits_{1\leq i<j\leq n}(j-i)=n(n^2-1)/6,
 $$
  the homogeneous dimension of  $T_1(n,\mathbb{R})$ equals to $n(n^2-1)/6.$  So, corollary 2 is valid for $T_1(n,\mathbb{R})$ with $Q=n(n^2-1)/6.$

\

\subsection{Finite-dimensional spaces over locally compact division rings}

\

Let $K$ be a  locally compact $\sigma$-compact division ring equipped with the norm $|\cdot|$ (e.g., $K=\mathbb{R}, \mathbb{Q}_p,$ or $\mathbb{H},$ the quaternion division ring). In the following we assume that the additive group $K^n$ is endowed with the invariant metric $\rho(x,y)=|x-y|_\infty:= \max_{1\leq i\leq n}|x_i-y_i|$ (here $x=(x_i),$ $y=(y_i)\in K^n$).

\textbf{Remark 5.} If $K$ is a field we have  \cite[Subsection VII.1.10, Corollary 1]{Burb}
 $$
 \mathrm{mod}(A(u))=\mathrm{mod}_K(\det(A(u))).
 $$

\textbf{Lemma 6.} \textit{The additive group $K^n$ endowed with the metric $\rho$ is a space of homogeneous type with respect to the Haar measure $\mu_n$ and }$C_{\mu_n}=2^n.$

Proof. First note that the additive group $K$ endowed with the metric $\rho_1(x,y)=|x-y|$ is a space of homogeneous type with respect to the Haar measure $\mu$ and $C_{\mu}=2.$ Indeed, since $(2e)B(0,r)=B(0,2r)$ ($e$ denotes the unit in $K$), we have for all $x\in K$ (see, e.g., \cite[Section VII.1, formula (32) and Definition 6]{Burb})
$$
\mu(B(x,2r))=\mu(B(0,2r))=\mathrm{mod}(2e)\mu(B(0,r))=2\mu(B(0,r))=2\mu(B(x,r)).
$$
Since $B(x,r)= \times_{i=1}^nB(x_i,r)$ where  $x=(x_1,\dots,x_n)\in K^n, r>0,$ lemma 6 follows.

 \textbf{Lemma 7.} \textit{For every family $A(u)$ of invertible $n\times n$ matrices with entries from $K$ the condition $(*)$ is valid with}
  $$
  k(u)=\max_{1\leq i\leq n}\sum\limits_{j=1}^n|a_{ij}(u)|
  $$
 \textit{ where $A(u)^{-1}=(a_{ij}(u))$ (in other words, $k(u)=\|A(u)^{-1}\|_\infty.$)
}

 Proof.  Indeed,   $A(u)^{-1}(B(x,r))=A(u)^{-1}x+A(u)^{-1}(B(0,r))$ and
$$
A(u)^{-1}(B(0,r))=\left\{\left(\sum\limits_{j=1}^na_{ij}(u)y_j\right)_{i=1}^n: y=(y_j)\in B(0,r)\right\}.
$$
Since
$$
\left|\left(\sum\limits_{j=1}^na_{ij}(u)y_j\right)_{i=1}^n\right|_\infty= \max_{1\leq i\leq n}\left|\sum\limits_{j=1}^na_{ij}(u)y_j\right|\leq  \max_{1\leq i\leq n}\sum\limits_{j=1}^n|a_{ij}(u)|r,
$$
we have $A(u)^{-1}(B(0,r))\subseteq B(0,k(u)r).$

\textbf{Definition 3.} We define the Hausdorff operator $\mathcal{H}$ on the additive group $G=K^n$ via the formula (1) where $A(u)$ is a family  of invertible $n\times n$ matrices with entries from $K$ and $\mu$ is replaced by $\mu_n.$

  %So, the Ces\`{a}ro operators in this case take the form
  %$$
%(\mathcal{C}_Af)(x)=\int\limits_{\{\mathrm{mod}_K(\det A(u))\geq 1\}} f(A(u)(x))\mathrm{mod}_K(\det A(u))d\mu(u).
%$$

Now theorem 1 along with remark 3 yield the next result.

 \textbf{  Corollary 3} (cf. \cite{LL}). \textit{The operator $\mathcal{H}$ is bounded on $H^1(K^n)$ provided $\Phi\in L^1_{k^n}(K^n)$ and}
$$
\|\mathcal{H}\|\leq 2^n\|\Phi\|_{L^1_{k^n}}.
$$

\textbf{Remark 6.} The  result by Ruan and Fan \cite[Theorem 1.3]{RF} shows that for $K=\mathbb{R}$ the above condition of boundedness for $\mathcal{H}$ can not be sharpened in general.

\textbf{Acknowledgments.} The author   thanks  E. Liflyand for posing the problem.

 This is a pre-print of the article \cite{Hausd}.

 Adolf R. Mirotin,

Department of Mathematics and   Programming  Technologies,
F. Skorina Gomel State University, 246019, Sovietskaya,
104, Gomel, Belarus

E-mail address: amirotin@yandex.ru
\end{document}